\documentclass{amsart}

\usepackage[margin=1in]{geometry}

\newtheorem{theorem}{Theorem}[section]
\newtheorem{lemma}[theorem]{Lemma}

\theoremstyle{definition}

\newtheorem{example}[theorem]{Example}

\theoremstyle{remark}

\numberwithin{equation}{section}

\usepackage{booktabs}
\usepackage{siunitx}
\usepackage{graphicx}
\usepackage{hyperref}
\usepackage{pdflscape}
\usepackage{xcolor}

\newtheorem{corollary}{Corollary}[theorem] 

\usepackage{listings}
\usepackage{xcolor}

\usepackage{amsmath}
\usepackage{amssymb}

\begin{document}

\title{On the Bounds of Stopping and Cycle Numbers in the Collatz Process}

\author{Daohang Sha}
\address{}
\curraddr{Philadelphia, PA 19104}
\email{dhsha@hotmail.com}
\thanks{}

\date{Oct 28, 2026}


\begin{abstract}
	We study the stopping-time and cycle structure of the normalized Collatz iteration using a binary sequence formalism and a recursive rooted sequence tree introduced by Winkler~\emph{(2017)}, we establish a bijection between admissible binary sequences and integers $m \equiv 3 \pmod{4}$, showing that every such integer arises uniquely in this construction. Within this framework, we derive sharp lower and upper bounds for the weighted digit-sum $\sum_{i=1}^{r} 3^{r-i} 2^{\mathbf{1}(i)-1}$, namely
	\[
	3^r - 2^r \leq \sum_{i=1}^{r} 3^{r-i} 2^{\mathbf{1}(i)-1} \leq r \cdot 3^{r-1},
	\]
	the upper bound being strictly tighter than the bound $s \cdot 3^r$ previously established by Luca~\emph{(2005)}. These bounds yield corresponding sharp bounds on any hypothetical cycle number $m_c$, and imply that the ratio $F_q(m)/m$ satisfies
	\[
	1 - \varepsilon_1 < \frac{F_q(m)}{m} < 1 + \varepsilon_2,
	\]
	with $F_q(m)/m \to 1$ as $r \to \infty$, while the ratio $3^r/2^s$ approaches but never reaches $1$ for any finite admissible sequence. We further show that the ratio of the maximum to minimum summation satisfies $\alpha(r) \sim r/3$ as $r \to \infty$, providing a theoretical basis for the empirically observed scaling $\alpha \approx \max(r)/3$. Together, these results imply that any nontrivial cycle number, if one exists, must grow without bound and satisfy increasingly stringent arithmetic constraints, consistent with all known computational and theoretical bounds on the Collatz conjecture.
\end{abstract}

\maketitle

\sloppy

\section{Introduction}
The $3n+1$ problem, also known as the Collatz conjecture, has been widely studied and remains one of the most intriguing unsolved problems in mathematics \cite{La11}, \cite{La12}. While several formulations of the Collatz function exist, we consider the following normalized form:
$$f(n)=\begin{cases}
	F_0(n)=\frac{n}{2}, & n = even \\
	F_1(n)=\frac{3n+1}{2}, & n = odd
\end{cases}$$

We denote the Collatz process as $F_q(n)=F_{101...0}(n)=F_0(...F_1(F_0(F_1(n))))$, where the binary sequence $q=101...0$, the total steps or the length of sequence is $s$ and the total number of odd steps is $r$. If $F_q(n)<n$, the minimum step $s$ is called the stopping time. If $F_q(m_c)=m_c$, $m_c$ is called a cycle number. The stopping time ($s$) and the cycle numbers ($m_c$) have been widely studied in \cite{Si10} - \cite{He23}. In \cite{Si10}, Simons and de Weger referred to a periodic sequence with $m$ local minima as an $m-cycle$. In the following discussion, we use the term "cycle" to denote any sequence that contains at least one such cycle. We also restrict our analysis to odd starting values $n \equiv 3 \pmod{4}$ satisfying $s \geq 4$, as those with $s \leq 2$ are readily confirmed by direct computation.

In this paper we obtain general lower and upper bounds for the stopping number and compute sharp upper bounds for both stopping and cycle numbers across all admissible sequences. Taken together, these bounds show that any nontrivial cycle would require an exact identity between powers of 2 and 3 in the expression for $F_q(n)$, an equality that cannot hold for any finite admissible sequence, implying that any nontrivial cycle, if one exists, must grow without bound and satisfy increasingly stringent arithmetic constraints.

\section{Collatz Sequence of Stopping}

Let $q$ be a binary sequence of length $s$, containing exactly $r$ entries equal to 1. We define the function $1(i)$ to denote the position of the $i^{\text{th}}$ occurrence of the digit 1 within the sequence $q$. For example, consider $n=423$, $F_q(n)=302$, the sequence $q = 1110110100$, which contains six ones. The corresponding positions of these ones are given by $1(1) = 1$, $1(2) = 2$, $1(3) = 3$, $1(4) = 5$, $1(5) = 6$, and $1(6) = 8$. Using this notation, the transformation $F_q(n)$ associated with the sequence $q$ in the Collatz process is given by
\begin{equation} \label{formula}
	F_q (n)=\dfrac{3^r n+\sum_{i=1}^r3^{r-i} 2^{1(i)-1}}{2^s},
\end{equation}

Here, $s$ is the total number of steps, $r$ is the number of applications of the function $F_1$, and $1(i)$ specifies the position of the $i^{\text{th}}$ occurrence of 1 in the sequence $q$.

\begin{lemma}\label{lem:length}
	Let \( r \in \mathbb{N} \) and define
	\[
	s = \left\lfloor 1 + r \cdot \log_2 3 \right\rfloor.
	\]
	Then there exist a non-negative integer \( i \in \mathbb{Z}_{\geq 0} \) and a positive integer \( k \in \mathbb{N} \) such that
	\[
	3^{r} \cdot (2i+1) - 1 = k \cdot 2^{s - r}.
	\]
\end{lemma}

\begin{proof}
	By the definition of \( s \), we have
	\[
	s = \left\lfloor 1 + r \cdot \log_2 3 \right\rfloor \quad \Rightarrow \quad s - 1 < r \log_2 3 < s.
	\]
	Exponentiating both sides with base 2 gives
	\[
	2^{s - 1} < 3^r < 2^s,
	\]
	which implies
	\[
	\frac{3^r}{2^{s - r}} \in \left(2^{r - 1}, 2^r\right).
	\]
	Now consider the congruence
	\[
	3^r \cdot t \equiv 1 \pmod{2^{s - r}}.
	\]
	Since \( 3^r \) is odd and \( \gcd(3^r, 2^{s - r}) = 1 \), the inverse of \( 3^r \) modulo \( 2^{s - r} \) exists. Therefore, there exists an odd integer \( t \equiv 2i + 1 \), for some \( i \in \mathbb{Z}_{\geq 0} \), such that
	\[
	3^r \cdot (2i + 1) \equiv 1 \pmod{2^{s - r}}.
	\]
	Rewriting this congruence yields
	\[
	3^r \cdot (2i + 1) - 1 = k \cdot 2^{s - r}
	\]
	for some \( k \in \mathbb{N} \), as required.
\end{proof}

\begin{lemma}\label{lem:min}
	Let \( r\geq2, i \in \mathbb{N} \), and for any starting number 
	\[
	m = 2^r \cdot (2i+1) - 1 \equiv 3\pmod{4}.
	\]
	Then the Collatz sequence of \( m \) begins with
	\[
	\underbrace{11\ldots1}_{r \text{ times}}\,\underbrace{00\ldots0}_{s - r \text{ times}}.
	\]
\end{lemma}

\begin{proof}
	Let us begin with the expression
	\[
	m = 2^r \cdot (2i+1) - 1.
	\]
	
	We claim that applying \( r \) iterations of the Collatz rule for odd integers, i.e., \( m \mapsto (3m + 1)/2 \), leads to
	\[
	F_{11...1}(m) = 3^{r} \cdot (2i+1) - 1.
	\]
	
	This is shown by induction:
	
	- Base case (step 1):
	\[
	\quad F_1(m) = \frac{3(2^r \cdot (2i+1) - 1) + 1}{2} = 3 \cdot 2^{r-1} (2i+1) - 1 = 2^{r-1} \cdot 3^{1} \cdot (2i+1) - 1.
	\]
	
	- Inductive step:
	Suppose that after \( r-1 \) steps,
	\[
	F_{11...1}(m) = 2 \cdot 3^{r - 1} \cdot (2i+1) - 1.
	\]
	Then, step $r$:
	\[
	F_{11...11}(m) = \frac{3[2 \cdot 3^{r - 1} \cdot (2i+1) - 1] + 1}{2} = \frac{2 \cdot 3^{r} \cdot (2i+1) - 2}{2} = 3^{r} \cdot (2i+1) - 1.
	\]
	
	So, after \( r \) such "odd steps" (each denoted by `1`), we reach \( 3^{r} (2i+1) - 1 \), which is even (since \( 3^{r} \) is odd, so the result is even). From this point, the Collatz process proceeds with divisions by 2 only, i.e., steps of the form \( m \mapsto m/2 \), which we denote by '0'.
	
	Hence, the step sequence begins with exactly \( r \) increasing steps ('1's), followed by decreasing steps ('0's) until the next odd number occurs.
	
	Therefore, together with the Lemma~\ref{lem:length}, the sequence begins with:
	\[
	\underbrace{11\ldots1}_{r \text{ times}}\,\underbrace{00\ldots0}_{s - r \text{ times}},
	\]
	completing the proof.
\end{proof}

\tiny
\begin{table}[ht]
	\caption{Specific Sequences of $q$ with $1(i)=i$ for $m \equiv 3\pmod{4}$ in Collatz process}
	\label{tab:seq}
\begin{tabular}{p{0.005\textwidth}p{0.005\textwidth}p{0.4\textwidth}p{0.11\textwidth}rp{0.125\textwidth}p{0.125\textwidth}}

		\toprule

		$r$ &  {$s$} & {Sequence($q$)} & {$2^r\cdot(2i+1)-1$} & {$m\pmod{2^s}$} & {$k=F_q(m)$} \\

		\midrule
		2	&	4	&	1100 & {$2^2 \cdot 1-1$} & 3 & 2\\
		3	&	5	&	11100 & {$2^3 \cdot 3-1$} & 23 & 20\\
		4	&	7	&	1111000 & {$2^4 \cdot 1-1$} & 15 & 10\\
		5	&	8	&	11111000 & {$2^5 \cdot 3-1$} & 95 & 91\\
		6	&	10	&	1111110000 & {$2^6 \cdot 9-1$} & 575 & 410\\
		7	&	12	&	111111100000 & {$2^7 \cdot 3-1$} & 383 & 205\\
		8	&	13	&	1111111100000 & {$2^8 \cdot 1-1$} & 255 & 205\\
		9	&	15	&	111111111000000	& {$2^9 \cdot 11-1$} & 5,631 & 3,383\\
		10	&	16	&	1111111111000000 & {$2^{10} \cdot 25-1$} & 25,599 & 23,066\\
		11	&	18	&	111111111110000000 & {$2^{11} \cdot 51-1$} & 104,447 & 70,582\\
		12	&	20	&	11111111111100000000 & {$2^{12} \cdot 17-1$} & 69,631 & 35,291\\
		13	&	21	&	111111111111100000000 & {$2^{13} \cdot 91-1$} & 745,471 & 566,732\\
		14	&	23	&	11111111111111000000000 & {$2^{14} \cdot 201-1$} & 3,293,183 & 1,877,689\\
		15	&	24	&	111111111111111000000000 & {$2^{15} \cdot 67-1$} & 2,195,455 & 1,877,689\\
		16	&	26	&	11111111111111110000000000 & {$2^{16} \cdot 193-1$} & 12,648,447 & 8,113,298\\
		17	&	27	&	111111111111111110000000000 & {$2^{17} \cdot 747-1$} & 97,910,783 & 94,206,740\\
		18	&	29	&	11111111111111111100000000000 & {$2^{18} \cdot 249-1$} & 65,273,855 & 47,103,370\\
		19	&	31	&	1111111111111111111000000000000	& {$2^{19} \cdot 83-1$}  & 43,515,903 & 23,551,685\\
		20	&	32	&	11111111111111111111000000000000 & {$2^{20} \cdot 1393-1$} & 1,460,666,367 & 1,185,813,152\\
		21	&	34	&	1111111111111111111110000000000000 & {$2^{21} \cdot 3195-1$} & 6,700,400,639 & 4,079,690,977\\
		22	&	35	&	11111111111111111111110000000000000 & {$2^{22} \cdot 1065-1$} & 4,466,933,759 & 4,079,690,977\\
		23	&	37	&	1111111111111111111111100000000000000 & {$2^{23} \cdot 8547-1$} & 71,697,432,575 & 49,111,434,902\\
		24	&	39	&	111111111111111111111111000000000000000 & {$2^{24} \cdot 2849-1$} & 47,798,288,383 & 24,555,717,451\\
		25	&	40	&	1111111111111111111111111000000000000000 & {$2^{25} \cdot 22795-1$} & 764,873,277,439 & 589,414,790,413\\
		26	&	42	&	111111111111111111111111110000000000000000 & {$2^{26} \cdot 18521-1$} & 1,242,923,270,143 & 718,351,699,928\\
		27	&	43	&	1111111111111111111111111110000000000000000 & {$2^{27} \cdot 28019-1$} & 3,760,646,520,831 & 3,260,217,528,257\\
		28	&	45	&	111111111111111111111111111100000000000000000 & {$2^{28} \cdot 31185-1$} & 8,371,159,695,359 & 5,442,907,506,622\\
		29	&	46	&	1111111111111111111111111111100000000000000000 & {$2^{29} \cdot 10395-1$} & 5,580,773,130,239 & 5,442,907,506,622\\
		30	&	48	&	111111111111111111111111111111000000000000000000 & {$2^{30} \cdot 3465-1$} & 3,720,515,420,159 & 2,721,453,753,311\\
		31	&	50	&	11111111111111111111111111111110000000000000000000 & {$2^{31} \cdot 263299-1$} & 565,430,297,034,751 & 310,197,425,018,629\\	
		
		\vdots	&	\vdots	&	\vdots & \vdots &\vdots &	\vdots 	\\	
		
		\bottomrule

	\end{tabular}
\end{table}
\normalsize

Table~\ref{tab:seq} includes several specific sequences for values of $r$ ranging from 2 to 31. This list can be extended indefinitely, as there is no upper bound on $r$. Therefore, there is no upper bound on $s$, too. For example, when $r=1000$, then $s=1585$,\\

$m=113181711197902800523896881066218160102709667045393956950389619432948098013957729384098\\
0894314243553691246663949249137551627877867250646145606647381546299358516204398946163838768783\\
5954105697727172101185270096853888059974345842688019492876298064702543347682952175404153595977\\
9733727568209962344081540465463827583787623071447065173716333421030611964595317030946266729588\\
1872408781238846687069393680025440093827021268188760756952836484710583975703867188025783792797\\
611927080534015$ (478 digits)\\

$F_q(m)=11027773766223939255599638243685383983123058040229118529895194635849895633361467411\\
1382402477750726619658150157035536038452913603129669634101848918205329307515073658919130294253\\
2781506097769243834069629788450644321665401676961180452393456153975938038355963675587842309211\\
9591171602528480255231243448888268364113078890677027695508662577537129844756989406360976065975\\
1963882667459405174174888685722485109676483971192668155890160044032425213449698721946740600735\\
2136396783848666645$ (478 digits)\\

$2i+1=105628373340941055460453010826461801731312915456661368786503022090718002344954512657\\4
11550899143065854323968254783269393108263739680866504502238016353334720601889819555845564641$ (177 digits)\\

\begin{theorem}\label{thm:min}
	The summation $ \sum_{i=1}^r 3^{r-i} \cdot 2^{1(i)-1} \geq 3^{r} - 2^{r},$ and equality holds when $1(i) = i$ for all $i$.
\end{theorem}

\begin{proof}
	To find the minimum value of
	$	\sum_{i=1}^r 3^{r-i} \cdot 2^{1(i)-1} $
	where $1(i)$ is the position of the $i^\text{th}$ 1 in a binary sequence $q$ of length $s$, we apply the result from the Rearrangement Inequality.	The minimum occurs when the positions $1(i)$ are as small as possible, i.e., $	1(i) = i \quad \text{for } i = 1, 2, \dots, r $
	That is, the 1s occur at the first $r$ positions in the sequence. Therefore, the minimum value is:
	$$	\sum_{i=1}^r 3^{r-i} \cdot 2^{i-1}	$$
    We can compute this closed-form as a geometric convolution:
	$$
	\sum_{i=1}^r 3^{r-i} \cdot 2^{i-1} = \sum_{j=0}^{r-1} 3^{j} \cdot 2^{r-1 - j} = 2^{r-1} \sum_{j=0}^{r-1} \left( \frac{3}{2} \right)^j
	$$
	This is a geometric series:	
	$$
	= 2^{r-1} \cdot \frac{ \left( \frac{3}{2} \right)^r - 1 }{ \frac{3}{2} - 1 } = 2^{r-1} \cdot \frac{ \left( \frac{3}{2} \right)^r - 1 }{ \frac{1}{2} } = 2^r \left[ \left( \frac{3}{2} \right)^r - 1 \right] = 3^{r} - 2^{r}
	$$
	Final Answer:
	$$
	\boxed {\sum_{i=1}^r 3^{r-i} \cdot 2^{1(i)-1} \geq 3^{r} - 2^{r} } 
	$$
	and equality holds when $1(i) = i$ for all $i$.
	
\end{proof}

\begin{theorem}[Maximizer]\label{thm:max}
	The summation 
	\[
	\sum_{i=1}^r 3^{\,r-i} \cdot 2^{\,1(i)-1} \le r\cdot 3^{r-1},
	\]
	attains its maximum value when 
	\[
	1(i) = \left\lfloor 1 + (i-1)\log_2 3 \right\rfloor, \quad i=1,\dots,r.
	\]
\end{theorem}

\begin{proof}	
	By Theorem~\ref{thm:exhaustive} (proved below), the tree in Figure 1 generates every admissible sequence exactly once; hence the maximum over tree sequences equals the maximum over all admissible sequences.
	
	Consider the modified rooted sequence tree introduced by Winkler \cite{Win17}, illustrated in  \autoref{fig:seqtree}. It is evident that the further to the right the $1$’s appear in a sequence, the larger resulting sum becomes comparing to the minimum summation when all $1$’s appear in the left. While top-line positions are the largest admissible positions consistent with the tree's insertion rule - see below for the details. Consequently, the sequences along the top line of the sequence tree yield the maximum values.
		
	\begin{figure}[ht]
		\centering
		\includegraphics[width=1.0\linewidth]{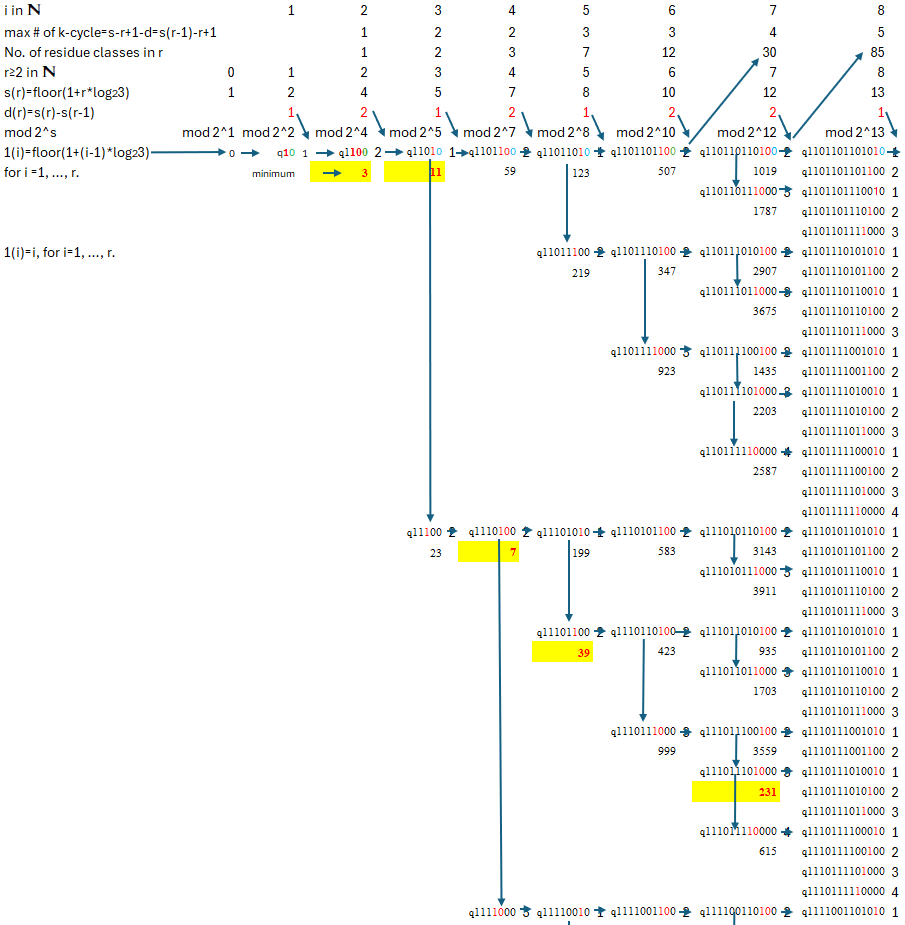}
		\caption{Rooted sequence tree}
		\label{fig:seqtree}
	\end{figure}

Here is how these sequences in the modified sequence tree are constructed. All sequences begin with $r=2, s=4$, so that the initial sequence is $q=1100$. This is also called rooted tree. Each sub-sequence is obtained by inserting either “10” or “1” immediately before each zero at the terminal zero-block in the previous sequence, according to the difference in their lengths. Specifically, define
\[
d(r) = s(r) - s(r-1) = 
\left\lfloor 1 + r \log_2 3 \right\rfloor - 
\left\lfloor 1 + (r-1)\log_2 3 \right\rfloor.
\]

This rule generates all admissible sub-sequences. And the total number of the admissible sequences in each $r$ equal to the summary of all number of zeros at the terminal zero-block in the previous sequence. Consequently, the positions of the (1)’s in the top-line sequence are given exactly by
\[
1(i) = \left\lfloor 1 + (i-1)\log_2 3 \right\rfloor, 
\quad i=1,\dots,r.
\]

Thus, for each $i$,
\[
1(i) - 1  = \big\lfloor 1+(i-1)\log_2 3\big\rfloor - 1
\le (i-1)\log_2 3,
\]
so
\[
2^{1(i)-1} \le 2^{\,(i-1)\log_2 3} = 3^{i-1}.
\]
Therefore every term of the sum satisfies
\[
3^{r-i}\cdot 2^{1(i)-1} \le 3^{r-i}\cdot 3^{i-1} = 3^{r-1}.
\]

Summing over $i = 1,\dots,r$,
\[
\sum_{i=1}^r 3^{r-i}2^{1(i)-1} \le \sum_{i=1}^r 3^{r-1} = r\cdot 3^{r-1}.
\]

This establishes the claim.
\end{proof}

	\textbf{Remark:} (Comparison to the literature)
	Recall the bound cited from Luca \cite{Lu05}:
	\[
	\sigma_q < \frac{s\cdot 3^r}{2^s}.
	\]
	Since $r < s$ for all admissible sequences with $r\ge 1$, the bound
	proven here,
	\[
	\sigma_q = \frac{\sum_{i=1}^r 3^{r-i}2^{1(i)-1}}{2^s} \le\ \frac{r\cdot 3^{r-1}}{2^s},
	\]
	is strictly tighter.

\begin{corollary}\label{thm:ratio}
	Let $\alpha(r)$ denote the ratio of the maximum summation bound 
	to the minimum summation, i.e.,
	\[
	\alpha(r) = \frac{r \cdot 3^{r-1}}{3^r - 2^r} 
	= \frac{r}{3\!\left(1 - (2/3)^r\right)}.
	\]
	Then $\alpha(r) > \dfrac{r}{3}$ for all finite $r \geq 2$, and
	\[
	\lim_{r \to \infty} \frac{\alpha(r)}{r/3} = 1.
	\]
	Consequently, $\alpha(r) \sim \dfrac{r}{3}$ as $r \to \infty$, 
	which validates the empirical observation $\alpha \approx \max(r)/3$.
\end{corollary}

\begin{proof}
	By Theorems~\ref{thm:min} and~\ref{thm:max}, the minimum and 
	maximum values of the summation $\sum_{i=1}^{r} 3^{r-i} 2^{\mathbf{1}(i)-1}$ 
	over all admissible sequences are $3^r - 2^r$ and $r \cdot 3^{r-1}$, 
	respectively. Their ratio is therefore
	\[
	\alpha(r) 
	= \frac{r \cdot 3^{r-1}}{3^r - 2^r} 
	= \frac{r}{3} \cdot \frac{1}{1-(2/3)^r}.
	\]
	Since $(2/3)^r \in (0,1)$ for all finite $r \geq 1$, the factor 
	$\bigl(1-(2/3)^r\bigr)^{-1} > 1$, which gives $\alpha(r) > r/3$.
	As $r \to \infty$, $(2/3)^r \to 0$, so 
	$\bigl(1-(2/3)^r\bigr)^{-1} \to 1$, and hence 
	\[
	\lim_{r \to \infty} \frac{\alpha(r)}{r/3} 
	= \lim_{r \to \infty} \frac{1}{1-(2/3)^r} = 1. \qedhere
	\]
\end{proof}

\textbf{Remark:}
	Corollary~\ref{thm:ratio} provides a theoretical justification for the 
	empirically observed relation $\alpha \approx \max(r)/3$: for large $r$, 
	the ratio of the maximum to minimum admissible summation is 
	asymptotically $r/3$, so across a dataset with largest odd-step count 
	$r_{\max}$, one expects $\alpha \approx r_{\max}/3$.  The strict 
	inequality $\alpha(r) > r/3$ for all finite $r$ reflects the fact that 
	the convergence is from above.

Table~\ref{tab:seq_max} shows the Collatz process of $r$ from 2 to 31, but with sequences with the maximum values of the summation. Again, the value of $r$ can be any big integers.

\tiny
\begin{table}[ht]
	\caption{Specific Sequences of $q$ with max value of the sum for $m \equiv 3\pmod{4}$ in Collatz process}
	\label{tab:seq_max}
		\begin{tabular}{p{0.005\textwidth}p{0.005\textwidth}p{0.4\textwidth}rp{0.13\textwidth}p{0.13\textwidth}}

			\toprule
			$r$ &  {$s$} & {Sequence($q$)} & {$m\pmod{2^s}$} &  {$F_q(m)<m$} \\
			\midrule
			2 &	4	& 1100     						& 3			&	2\\
			3 &	5	& 11010							& 11		&	10\\
			4 &	7	& 1101100						& 59		&	38\\
			5 &	8	& 11011010      				& 123		&	118\\
			6 &	10	& 1101101100	  				& 507		&	362\\
			7 &	12	& 110110110100  				& 1,019 	&	545\\
			8 &	13	& 1101101101010		 			& 3,067		&	2,458\\
			9 &	15	& 110110110101100 				& 7,163  	&	4,304\\
			10 & 16	& 1101101101011010 				& 56,315	&	50,743\\
			11 & 18	& 110110110101101100 	    	& 23,547	&	15,914\\
			12 & 20	& 11011011010110110100    		& 416,763  	&	211,226\\
			13 & 21	& 110110110101101101010   		& 941,051 	&	715,420\\
			14 & 23	& 11011011010110110101100  		& 8,281,083	&	4,721,663\\
			15 & 24	& 110110110101101101011010		& 4,086,779 &	3,495,268\\
			16 & 26	& 11011011010110110101101100	& 62,807,035 &	40,287,332\\
			17 & 27	& 110110110101101101011011010	& 96,361,467 &  92,716,039\\
			18 & 29	& 11011011010110110101101101100	& 431,905,787 & 311,674,835\\
			19 & 31	& 1101101101011011010110110110100 & 163,470,331 & 88,473,443\\
			20 & 32	& 11011011010110110101101101101010 & 1,237,212,155 & 1,004,406,265\\
			21 & 34	& 1101101101011011010110110110101100 & 7,679,663,099 & 4,675,937,150\\
			22 & 35	& 11011011010110110101101101101011010 &	33,449,466,875 & 30,549,700,432\\
			23 & 37	& 1101101101011011010110110110101101100	& 119,348,812,795 & 81,751,762,091\\
			24 & 39	& 110110110101101101011011011010110110100 &	462,946,196,475 & 237,832,281,869\\
			25 & 40	& 1101101101011011010110110110101101101010 & 188,068,289,531 & 144,926,270,443\\
			26 & 42	& 110110110101101101011011011010110110101100 & 2,936,847,358,971 & 1,697,360,845,538\\
			27 & 43	& 1101101101011011010110110110101101101011010 &	5,135,870,614,523 & 4,452,440,639,554\\
			28 & 45	& 110110110101101101011011011010110110101101100 & 18,330,010,147,835 & 11,918,127,650,276\\
			29 & 46	& 1101101101011011010110110110101101101011011010 & 35,922,196,192,251 & 35,034,785,816,635\\
			30 & 48	& 110110110101101101011011011010110110101101101100 & 211,844,056,636,411 & 154,958,046,921,632\\
			31 & 50	& 11011011010110110101101101101011011010110110110100 & 915,531,498,413,051 & 502,264,407,868,691\\
																
			\vdots	&	\vdots	&	\vdots & \vdots &\vdots 	\\	
			\bottomrule
		\end{tabular}
	\end{table}
\normalsize

\begin{lemma}\label{lem:pow2}
	For any integer $s \ge 1$,
	\[
	\sum_{i=0}^{s-1} 2^{\,i} = 2^{\,s} - 1.
	\]
\end{lemma}

\begin{proof}
	Consider the finite geometric series with ratio $2$:
	\[
	\sum_{i=0}^{s-1} 2^{\,i} = \frac{2^{\,s} - 1}{2 - 1} = 2^{\,s} - 1.
	\]
\end{proof}

\begin{lemma}\label{lem:3mod4}

Let $m \equiv 3 \pmod{4}$. Then $m$ admits a unique representation of the form
$$ m = 3 + \sum_{n_i\ge2} 2^{n_i}, \qquad n_i \ge 2, $$
where the exponents $n_i$ are distinct. Conversely, every finite sum of the form
$$ 3 + \sum_{n_i\ge2} 2^{n_i}, \qquad n_i \ge 2,$$
yields an integer $m \equiv 3 \pmod{4}$.

\end{lemma}

\begin{proof}

Since $2^{n} \equiv 0 \pmod{4}$ for all $n \ge 2$, any integer of the form $3 + \sum_{n_i} 2^{n_i}$ is congruent to $3\pmod{4}$. Conversely, if $m \equiv 3 \pmod{4}$, then $m-3$ is divisible by 4 and hence has a binary expansion involving only powers $2^{n}$ with $n \ge 2$. Uniqueness follows from the uniqueness of the binary expansion.

\end{proof}

\begin{example}\label{exp:mod4}
For instance,
\small
\[
\begin{array}{c|l|l|l}
	m=4i-1 \equiv 3\pmod{4} & m=\sum_{n_i \ge 2} 2^{n_i} + 3 & 2^{s}y + m & 3^{r}y + F_{q}(m) \\ \hline
	3  & 3 & 2^{4}\cdot0 + 3 & 3^{2}\cdot0 + 2\\
	7  & 2^{2} + 3 & 2^{7}\cdot0 + 7 & 3^{4}\cdot0 +5\\
	11 & 2^{3} + 3 & 2^{5}\cdot0 + 11 & 3^{3}\cdot0 +10\\
	15 & 2^{3} + 2^{2} + 3 & 2^{7}\cdot0 + 15 & 3^{4}\cdot0 +10\\
	19 & 2^{4} + 3 & 2^{4}\cdot1 + 3 & 3^{2}\cdot1 +11\\
	23 & 2^{4} + 2^{2} + 3 & 2^{5}\cdot0 + 23 & 3^{3}\cdot0 + 20\\
	27 & 2^{4} + 2^{3} + 3 & 2^{59}\cdot0 + 27 & 3^{37}\cdot0 + 23\\
	31 & 2^{4} + 2^{3} + 2^{2} + 3 & 2^{56}\cdot0 + 31 & 3^{35}\cdot0 + 23\\
	35 & 2^{5} + 3 & 2^{4}\cdot2 + 3 & 3^{2}\cdot2 + 2\\
	39 & 2^{5} + 2^{2} + 3 & 2^{8}\cdot0 + 39 & 3^{5}\cdot0 + 38\\
	43 & 2^{5} + 2^{3} + 3 & 2^{5}\cdot1 + 11 & 3^{3}\cdot1 + 10\\
	47 & 2^{5} + 2^{3} + 2^{2} + 3 & 2^{54}\cdot0 + 47 & 3^{34}\cdot0 + 46\\
	51 & 2^{5} + 2^{4} + 3 & 2^{4}\cdot3 + 3 & 3^{2}\cdot3 + 2\\
	55 & 2^{5} + 2^{4} + 2^{2} + 3  & 2^{5}\cdot1 + 23 & 3^{3}\cdot1 + 20\\
	59 & 2^{5} + 2^{4} + 2^{3} + 3 & 2^{7}\cdot0 + 59 & 3^{4}\cdot0 + 38\\
	63 & 2^{5} + 2^{4} + 2^{3} + 2^{2} + 3 & 2^{54}\cdot0 + 63 & 3^{34}\cdot0 + 61\\
	67 & 2^{6} + 3 & 2^{4}\cdot4 + 3 & 3^{2}\cdot4 + 2\\
	71 & 2^{6} + 2^{2} + 3 & 2^{51}\cdot0 + 71 & 3^{32}\cdot0 + 61\\
	75 & 2^{6} + 2^{3} + 3 & 2^{5}\cdot2 + 11 & 3^{3}\cdot2 + 10\\
	79 & 2^{6} + 2^{3} + 2^{2} + 3 & 2^{8}\cdot0 + 79 & 3^{5}\cdot0 + 76\\
	83 & 2^{6} + 2^{4} + 3 & 2^{4}\cdot5 + 3 & 3^{2}\cdot5 + 2\\
	87 & 2^{6} + 2^{4} + 2^{2} + 3  & 2^{5}\cdot2 + 23 & 3^{3}\cdot2 + 20\\
	91 & 2^{6} + 2^{4} + 2^{3} + 3 & 2^{45}\cdot0 + 91 & 3^{28}\cdot0 + 61\\
	95 & 2^{6} + 2^{4} + 2^{3} + 2^{2} + 3 & 2^{8}\cdot0 + 95 & 3^{5}\cdot0 + 91\\
	99 & 2^{6} + 2^{5} + 3 & 2^{4}\cdot6 + 3 & 3^{2}\cdot6 + 2\\
   103 & 2^{6} + 2^{5} + 2^{2} + 3 & 2^{42}\cdot0 + 103 & 3^{26}\cdot0 + 61\\
   107 & 2^{6} + 2^{5} + 2^{3} + 3 & 2^{5}\cdot3+11 & 3^{3}\cdot3 + 10\\
   111 & 2^{6} + 2^{5} + 2^{3} + 2^{2} + 3 & 2^{31}\cdot0+111 & 3^{19}\cdot0 + 61\\
   115 & 2^{6} + 2^{5} + 2^{4} + 3 & 2^{4}\cdot7 + 3 & 3^{2}\cdot7 + 2\\
   119 & 2^{6} + 2^{5} + 2^{4} +2^{2} + 3 & 2^{5}\cdot3+23 & 3^{3}\cdot3 + 20\\
   123 & 2^{6} + 2^{5} + 2^{4} +2^{3} + 3 & 2^{8}\cdot0+123 & 3^{5}\cdot0 + 118\\
   127 & 2^{6} + 2^{5} + 2^{4} +2^{3} + 2^{2} + 3 & 2^{15}\cdot0+127 & 3^{9}\cdot0 +77 \\	
   131 & 2^{7}+3 & 2^{4}\cdot8+3 & 3^{2}\cdot8 + 2\\
   135 & 2^{7}+2^{2}+3 & 2^{7}\cdot0+135 & 3^{4}\cdot0 + 86\\
   139 & 2^{7}+2^{3}+3 & 2^{5}\cdot4+11 & 3^{3}\cdot4 + 10\\
   143 & 2^{7}+2^{3}+2^{2}+3 & 2^{7}\cdot1+15 & 3^{4}\cdot1 + 10\\
   147 & 2^{7}+2^{4}+3 & 2^{4}\cdot9+3 & 3^{2}\cdot9 + 2\\
   151 & 2^{7}+2^{4}+2^{2}+3 & 2^{5}\cdot4+23 & 3^{3}\cdot4 + 20\\   
   155 & 2^{7}+2^{4}+2^{3}+3 & 2^{40}\cdot0+155 & 3^{25}\cdot0 + 122\\
   159 & 2^{7}+2^{4}+2^{3}+2^{2}+3 & 2^{21}\cdot0+159 & 3^{13}\cdot0 + 122\\
   163 & 2^{7}+2^{5}+3 & 2^{4}\cdot10+3 & 3^{2}\cdot10 + 2\\
   167 & 2^{7}+2^{5}+2^{2}+3 & 2^{29}\cdot0+167 & 3^{18}\cdot0 + 122\\
   171 & 2^{7}+2^{5}+2^{3}+3 & 2^{5}\cdot5+11 & 3^{3}\cdot5 + 10\\
   175 & 2^{7}+2^{5}+2^{3}+2^{2}+3 & 2^{8}\cdot0+175 & 3^{5}\cdot0 + 167\\
   179 & 2^{7}+2^{5}+2^{4}+3 & 2^{4}\cdot11+3 & 3^{2}\cdot11+ 2\\
   183 & 2^{7}+2^{5}+2^{4}+2^{2}+3 & 2^{5}\cdot5+23 & 3^{3}\cdot5 + 20\\
   187 & 2^{7}+2^{5}+2^{4}+2^{3}+3 & 2^{7}\cdot0+187 & 3^{4}\cdot0 + 119\\
   191 & 2^{7}+2^{5}+2^{4}+2^{3}+2^{2}+3 & 2^{13}\cdot0+191 & 3^{8}\cdot0 + 154\\
   195 & 2^{7}+2^{6}+3 & 2^{4}\cdot12+3 & 3^{2}\cdot12 + 2\\
   199 & 2^{7}+2^{6}+2^{2}+3 & 2^{8}\cdot0+199 & 3^{5}\cdot0 + 190\\
   203 & 2^{7}+2^{6}+2^{3}+3 & 2^{5}\cdot6+11 & 3^{3}\cdot6 + 10\\
   207 & 2^{7}+2^{6}+2^{3}+2^{2}+3 & 2^{13}\cdot0+207 & 3^{8}\cdot0 + 167\\
   211 & 2^{7}+2^{6}+2^{4}+3 & 2^{4}\cdot13+3 & 3^{2}\cdot13 + 2\\
   215 & 2^{7}+2^{6}+2^{4}+2^{2}+3 & 2^{5}\cdot6+23 & 3^{3}\cdot6 + 20\\
   219 & 2^{7}+2^{6}+2^{4}+2^{3}+3 & 2^{8}\cdot0+219 & 3^{5}\cdot0 + 209\\
   223 & 2^{7}+2^{6}+2^{4}+2^{3}+2^{2}+3 & 2^{31}\cdot0+223 & 3^{19}\cdot0 + 122\\
   227 & 2^{7}+2^{6}+2^{5}+3 & 2^{4}\cdot14+3 & 3^{2}\cdot14 + 2\\
   231 & 2^{7}+2^{6}+2^{5}+2^{2}+3 & 2^{12}\cdot0+231 & 3^{7}\cdot0 + 124\\
   235 & 2^{7}+2^{6}+2^{5}+2^{3}+3 & 2^{5}\cdot7+11 & 3^{3}\cdot7 + 10\\
   239 & 2^{7}+2^{6}+2^{5}+2^{3}+2^{2}+3 & 2^{20}\cdot0+239 & 3^{12}\cdot0 + 122\\
	\vdots & \vdots & \vdots  & \vdots
\end{array}
\]
\normalsize
\end{example}

\begin{lemma}[Uniqueness of sequences and starting values]\label{lem:seqtree}
	The binary sequences $q$ displayed in Figure~\ref{fig:seqtree} are pairwise distinct. 
	Moreover, the construction defines an injective mapping
	\[
	q \longmapsto m,
	\]
	where each sequence $q$ corresponds to a unique integer $m \equiv 3 \pmod{4}$
	satisfying $F_q(m) < m$.
\end{lemma}

\begin{proof}
	The sequences in Figure~\ref{fig:seqtree} are generated recursively using a
	deterministic insertion rule: at each stage $r$, a new sequence is obtained from a
	unique parent sequence by inserting either the substring $1$ or $10$ at a prescribed
	position immediately preceding the terminal block of zeros. Since both the parent
	sequence and the insertion position are uniquely determined, distinct construction
	paths necessarily produce distinct binary sequences. Hence all sequences $q$ in
	Figure~\ref{fig:seqtree} are pairwise distinct.
	
	For a fixed sequence $q$ of length $s$, equation~\eqref{formula} uniquely determines
	the residue class of $m \pmod{2^s}$. By Lemma~\ref{lem:3mod4}, every integer
	$m \equiv 3 \pmod{4}$ admits a unique binary expansion of the form
	\[
	m = 3 + \sum_{n_i \ge 2} 2^{n_i}.
	\]
	The exponents $n_i$ correspond uniquely to the positions of the $1$'s in the
	sequence $q$. Consequently, distinct sequences $q$ give rise to distinct integers
	$m$. Finally, direct comparison of the iterates of $m$ with its binary expansion shows
	that the constructed integer satisfies $F_q(m) < m$. This establishes injectivity
	of the mapping $q \mapsto m$.
\end{proof}

\begin{example}	
For instance, let $i$ is any integer and $4i-1 \equiv 3\pmod{4}$, then the first two iterates are
\[
F_{1}(4i-1) = 3\cdot 2i - 1, 
\qquad 
F_{11}(4i-1) = 3^2 i - 1.
\]

As an illustrative example, take 
\[
r = 11, \quad s = 18, \quad q = 111111111010100000.
\]
Then the iterates proceed as
\[
\begin{aligned}
	F_{111}(4i-1) &= \frac{3(3^2(2j)-1)+1}{2} = 3^3 j - 1,\\
	F_{1111}(4i-1) &= \frac{3(3^3(2k)-1)+1}{2} = 3^4 k - 1,\\
	F_{11111}(4i-1) &= \frac{3(3^4(2\ell)-1)+1}{2} = 3^5 \ell - 1,\\
	F_{111111}(4i-1) &= \frac{3(3^5(2m)-1)+1}{2} = 3^6 m - 1,\\
	F_{1111111}(4i-1) &= \frac{3(3^6(2n)-1)+1}{2} = 3^7 n - 1,\\
	F_{11111111}(4i-1) &= \frac{3(3^7(2o)-1)+1}{2} = 3^8 o - 1,\\
	F_{111111111}(4i-1) &= \frac{3(3^8(2p)-1)+1}{2} = 3^9 p - 1,\\
	F_{1111111110}(4i-1) &= \frac{3^9(2q+1)-1}{2} = 3^9 q + 9841,\\
	F_{11111111101}(4i-1) &= \frac{3(3^9(2r)+9841)+1}{2} = 3^{10} r + 14762,\\
	F_{111111111010}(4i-1) &= \frac{3^{10}(2s)+14762}{2} = 3^{10} s + 7381,\\
	F_{1111111110101}(4i-1) &= \frac{3(3^{10}(2t)+7381)+1}{2} = 3^{11} t + 11072,\\
	F_{11111111101010}(4i-1) &= \frac{3^{11}(2u)+11072}{2} = 3^{11} u + 5536,\\
\end{aligned}	
\]
\[
\begin{aligned}
	F_{111111111010100}(4i-1) &= \frac{3^{11}(2v)+5536}{2} = 3^{11} v + 2768,\\
	F_{1111111110101000}(4i-1) &= \frac{3^{11}(2w)+2768}{2} = 3^{11} w + 1384,\\
	F_{11111111101010000}(4i-1) &= \frac{3^{11}(2x)+1384}{2} = 3^{11} x + 692,\\
	F_{111111111010100000}(4i-1) &= \frac{3^{11}(2y)+692}{2} = 3^{11} y + 346 = 3^r y + F_q(m).
\end{aligned}
\]

The corresponding expansion of $4i-1$ in powers of two is
\[
\begin{aligned}
	4i-1 &= 2^2 (2j) - 1 = 2^2 2^2 k - 1 = 2^2 2^3 \ell - 1 = 2^2 2^4 m - 1 = 2^2 2^5 n - 1 = 2^2 2^6 o - 1\\
	&= 2^2 2^7 p - 1 = 2^2 2^7 (2q+1) - 1 = 2^2 2^8 q + 2^2 2^7 - 1 = \hdots = 2^{18} y + 511\\
	&= 2^s y + \sum_{n_i \ge 2} 2^{n_i} + 3 = 2^s y + m.\\
\end{aligned}
\]
\end{example}

Note that \(q, r, s\) appearing in the iteration steps are simply index variables; they are not the number of ones (\(r\)) or the sequence length (\(s\)) associated with the specific string \(q = 111111111010100000\). Their roles are distinct in context without confusion.

\begin{theorem}[Exhaustivity for $m \equiv 3 \pmod{4}$]\label{thm:exhaustive}
	The construction in Figure~\ref{fig:seqtree} generates each integer
	$m \equiv 3 \pmod{4}$ exactly once. Equivalently, the correspondence
	\[
	q \longleftrightarrow m
	\]
	defined by Figure~\ref{fig:seqtree} is a bijection between admissible sequences $q$
	and integers $m \equiv 3 \pmod{4}$.
\end{theorem}

\begin{proof}
	We prove the statement by showing that the constructions
	$q \mapsto m$ and $m \mapsto q$ are well-defined and inverse to each other.
	
	\medskip
	\noindent\textbf{Forward direction ($q \mapsto m$).}
	Starting from the root $q = 11$, we have
	\[
	F_{11}(4i - 1) = 3^{2} i - 1.
	\]
	Any admissible sequence $q = 11\cdots$ corresponds to a path in the rooted
	sequence tree shown in Figure~\ref{fig:seqtree}. Along this path, the Collatz
	process, $F_q(4i - 1)=3^ry+F_q(m)$, yields successive decompositions
	\[
	i = 2j \text{ or } 2j + 1,\quad
	j = 2k \text{ or } 2k + 1,\quad
	k = 2l \text{ or } 2l + 1,\ \ldots,\ 
	x = 2y \text{ or } 2y + 1.
	\]
	Consequently, we may write
	\[
	4i - 1 = 2^{s} y + \sum_{n_i \ge 2} 2^{n_i} + 3
	= 2^{s} y + m,
	\]
	which implies
	\[
	m \equiv 3 \pmod{4}.
	\]
	
	\medskip
	\noindent\textbf{Backward direction ($m \mapsto q$).}
	Let $m$ be any integer satisfying $m \equiv 3 \pmod{4}$. By
	Lemma~\ref{lem:3mod4}, $m$ admits a representation of the form
	\[
	m = \sum_{n_i \ge 2} 2^{n_i} + 3.
	\]
	Writing $4i - 1 = 2^ry + m$ and using the powers $2^{n_i}$ in the above decomposition,
	we recover successive relations
	\[
	i = 2j \text{ or } 2j + 1,\quad
	j = 2k \text{ or } 2k + 1,\quad
	k = 2l \text{ or } 2l + 1,\ \ldots,
	\]
	which uniquely determine an admissible sequence
	\[
	q = 11\cdots
	\]
	can be obtained via the inverse Collatz process and the mapping $F_q(4i - 1)=3^ry+F_q(m)$.
	
	\medskip
	Since the mappings $q \mapsto m$ and $m \mapsto q$ are inverses of each other,
	the correspondence between admissible sequences $q$ and integers
	$m \equiv 3 \pmod{4}$ is bijective.
\end{proof}

\begin{theorem}\label{thm:loop}
	For any integer $i$, we have
	\[
	4i - 1 \equiv 3 \pmod{4}, \quad \text{and} \quad F_q(4i - 1) \neq 4i - 1
	\]
	for all cases where $y \neq 0$.
\end{theorem}

\begin{proof}
	By definition, we have
	\[
	4i - 1 = 2^{s} y + m,
	\]
	and
	\[
	F_q(4i - 1) = 3^{r} y + F_q(m),
	\]
	where $y$ is an integer.
	
	Suppose, for the sake of contradiction, that
	\[
	4i - 1 = F_q(4i - 1).
	\]
	Assuming $m = F_q(m)$, it follows by substitution that
	\[
	2^{s} y = 3^{r} y.
	\]
	Since $y \neq 0$ and $r, s \geq 1$ by definition, dividing both sides by $y$ yields
	\[
	2^{s} = 3^{r}.
	\]
	This is impossible, as no positive integer power of $2$ can equal a power of $3$. Therefore, $F_q(4i - 1) \neq 4i - 1$, which completes the proof.
\end{proof}

\textbf{Remark:} This implies that for any sequence where $y \neq 0$, the existence of a loop is impossible. Consequently, a fixed point satisfying $F_{q}(m)=m$ can only potentially exist when $y=0$.

Let $T$ denote the entire set of elements represented in the table of Example~\ref{exp:mod4}. We partition $T$ into two disjoint sets, $G_0$ and $G_{\neq 0}$, based on the value of $y$:
$$
	G_0 = \left\{ \left( 2^s \cdot 0 + m, \, 3^r \cdot 0 + F_q(m) \right) \;\middle|\; y = 0 \right\} $$ 
$$
	G_{\neq 0} = \left\{ \left( 2^s y + m, \, 3^r y + F_q(m) \right) \;\middle|\; y \neq 0 \right\}
$$
where $s, r, m,$ and $y$ satisfy the conditions specified by the table columns. 

By construction, these two sets partition the table; they are mutually disjoint:
$$	G_0 \cap G_{\neq 0} = \emptyset $$
and their union constitutes the complete table:
$$	G_0 \cup G_{\neq 0} = T $$

Theorem~\ref{thm:loop} establishes that a fixed point $F_{q}(m)=m$ is impossible within $G_{\neq 0}$. Crucially, this impossibility extends to $G_0$ as well. Every base element $2^s\cdot0+m, 3^r\cdot0+F_q(m)$ in $G_0$ generates a family of elements $2^sy+m, 3^ry+F_q(m)$ in $G_{\neq 0}$ via the exact same parameters $s, r,$ and $y \neq 0$. If $F_{q}(m)=m$ were to hold for an element in $G_0$, then the corresponding relation $F_q(m)=m$ would simultaneously hold for its generated elements in $G_{\neq 0}$. This directly contradicts Theorem~\ref{thm:loop}. Therefore, fixed points cannot exist in either $G_0$ or $G_{\neq 0}$.

\section{Bounds of Stopping Number and Cycle Number}

Suppose that there exists a cycle in the sequence, i.e. $$F_q (m_c )=m_c.$$ Using equation \eqref{formula}, we have $$\dfrac{3^rm_c+\sum_{i=1}^r3^{r-i} 2^{1(i)-1}}{2^s}=m_c.$$
Thus, 
\begin{equation*}
	m_c=\dfrac{\sum_{i=1}^r3^{r-i} 2^{1(i)-1}}{2^s-3^r}.
\end{equation*}

\begin{theorem}[Sharpened cycle-number bound]\label{thm:mc_bound}
	If a nontrivial cycle number
	$m_c$ exists for parameters $(r,s)$, then
	\[
	\frac{3^r - 2^r}{2^s - 3^r} \le m_c \le \frac{r\cdot 3^{r-1}}{2^s - 3^r}.
	\]
\end{theorem}

\begin{proof}
Using \autoref{thm:min} and \autoref{thm:max}, we immediately obtain

\begin{align*}
m_c=\dfrac{\sum_{i=1}^r3^{r-i} 2^{1(i)-1}}{2^s-3^r}\geq\min_{q} \dfrac{\sum_{i=1}^r3^{r-i} 2^{1(i)-1}}{2^s-3^r}
=\dfrac{3^{r}-2^{r}}{2^s-3^r}, 
\end{align*}

and 

\begin{align*}
	m_c=\dfrac{\sum_{i=1}^r3^{r-i} 2^{1(i)-1}}{2^s-3^r}\leq\max_{q} \dfrac{\sum_{i=1}^r3^{r-i} 2^{1(i)-1}}{2^s-3^r}
	=\dfrac{r \cdot 3^{r-1}}{2^s-3^r}, 
\end{align*}

as claimed.

\end{proof}

\textbf{Remark:} Under the condition of Lemma~\ref{lem:length}, i.e., $ s = \left\lfloor 1 + r \cdot \log_2 3 \right\rfloor$, thus  $3^r<2^s$ holds, which makes the denominator $2^s-3^r$ positive.

For example, the value of expression
$$
\frac{3^r - 2^r}{2^s - 3^r}
\quad \text{for } r = 431{,}166{,}034{,}846{,}567 \text{ and } s = 683{,}381{,}996{,}816{,}440
$$
is approximately: \boxed{7.50 \times 10^{14}}. The values of $r$ and $s$ used here correspond to the maximum values listed in the reference \cite{El93}. Based on the results from \cite{Lu05}, the upper bound for the number of non-trivial cycles can be expressed as
 $$
\boxed { m_c < \dfrac{s \cdot 3^{r}}{2^s - 3^r}. }
 $$
 Substituting $r = 431{,}166{,}034{,}846{,}567$ and $s = 683{,}381{,}996{,}816{,}440$, this upper bound evaluates to approximately  \boxed{5.12 \times 10^{29}}.

\begin{figure}[tb]
	\centering
	\includegraphics[scale=0.55]{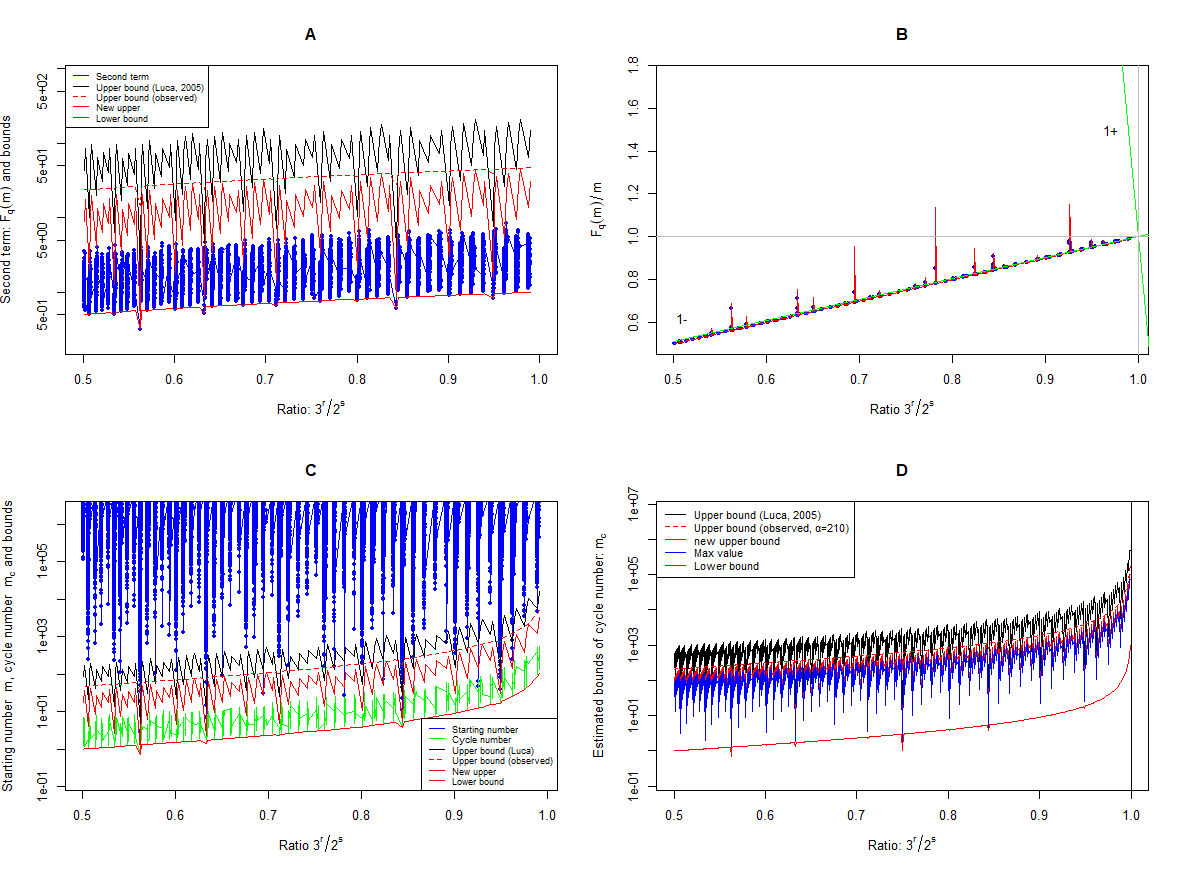}
	\caption{Collatz process results for integers $m \equiv 3\pmod{4}$, ranging from 3 to 16,000,000. A. Second term $\sigma_q$ (blue dots) and its bounds in green and red versus $3^r/2^s$; B. Ratio $F_q (m)/m$ (blue dots) versus $3^r/2^s$; C. Starting number $m$ (blue dots), cycle number $m_c$ in green and its bounds in red and black versus $3^r/2^s$; D. Estimated bounds for the cycle number $m_c$ versus $3^r/2^s$ using extended $r$.}
	\label{fig3}
\end{figure}

Rewrite the stopping number as,
	$$F_q (m)=\dfrac{3^r m+\sum_{i=1}^r3^{r-i} 2^{1(i)-1}}{2^s}=\dfrac{3^r}{2^s}m +\sigma_q,$$
where $$\sigma_q=\dfrac{\sum_{i=1}^r3^{r-i} 2^{1(i)-1}}{2^s},$$
it is easy to see that the stopping number $F_q (m)$ consists of two terms, the first term includes starting number $m$ and is unrelated to the sequence $q$, while the second term doesn’t include starting number $m$, but it is only related to the sequence. The value of $F_q (m)$ is dominated by the first term. The second term determines the position of $F_q (m)$ between $m$ and $m/2$. Similar to the the cycle number, the lower bound of the second term $\sigma_q$ is determined by
$$\boxed{\sigma_q\geq\dfrac{3^{r}-2^{r}}{2^s}, }$$
while the upper bound based on the results of \cite{Lu05} is
$$\boxed{\sigma_q<\dfrac{s3^{r}}{2^s}, }$$
and from the observation data that are showing in \autoref{fig3}. Here, the starting numbers $m$ for \autoref{fig3}A, B, and C are from the positive integer of $3\pmod{4}$ start from 3 to 16,000,000 whose stopping time $s\geq4$. We can see that 
$$\boxed{\sigma_q<\alpha\dfrac{3^{r}-2^{r}}{2^s},}$$
where $\alpha$ is a constant dependent on the maximum value of $r$ in all sequences.

\begin{figure}[tb]
	\centering
	\includegraphics[scale=0.55]{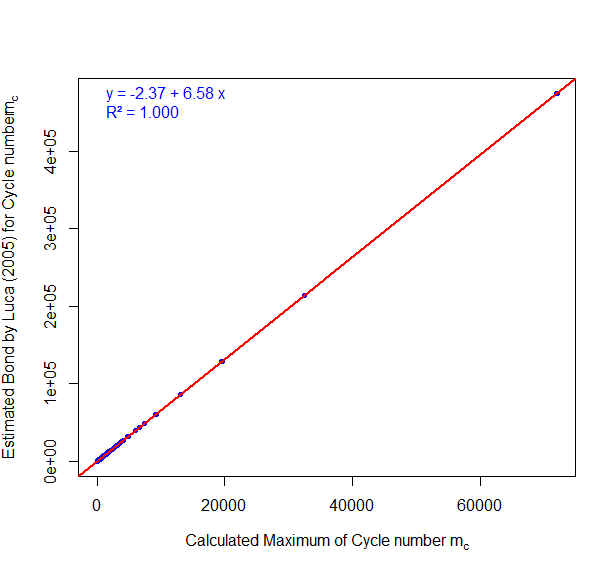}
	\caption{Regression of Estimated bound for the cycle number $m_c$ versus Calculated max values with $max(r)=630$.}
	\label{reg_mc}
\end{figure}

From the computing results, we observed that \boxed{\alpha=max(r)/3}. Here $\alpha$=10 for  \autoref{fig3}A,B and C; $\alpha$=210 for \autoref{fig3}D with max of $r=630$. We also note that the estimated bond by Luca (2005) for cycle numbers in \cite{Lu05} has the same pattern as our calculated max of cycle numbers that plotted in the \autoref{fig3}D. The relationship between them can be regressed as  \boxed{y = 6.58 x -2.37}, which is plotted in the \autoref{reg_mc}.

\autoref{fig3}A shows the second term of $F_q (m)$ (i.e. $\sigma_q$) and its bounds versus the ratio of $3^r$ and $2^s$. The \autoref{fig3}B shows that the non-trivial cycle number would be where $\dfrac{3^{r}}{2^s}$ is close to 1 if it existed. \autoref{fig3}B shows the ratio $F_q (m)/m$ versus ratio of $3^r$ and $2^s$. \autoref{fig3}C shows the starting number $m$ (blue dots) and the cycle number $m_c$ (green line) and its bounds (red lines). \autoref{fig3}D is showing the bounds of the cycle number with more values of $3^r/2^s$. The data in this figure was generated for $r=1$ to 630 (645 is the maximum number that R software can go with), then we apply $ s = \left\lfloor 1 + r \cdot \log_2 3 \right\rfloor$. The value of $s$ corresponding to $r$ is unique because of this restriction. Then we use the pair of $(r,s)$ to calculate the cycle number bound. We found that the minimum cycle number $m_c$ is greater than 978 when $(r,s)=(306,485)$ and $3^r/2^s\approx0.9989783$. It is not hard to see that the ratio $F_q (m)/m$ tends to be 1 as ratio $3^r/2^s$ goes to 1. Therefore, if there exists a cycle, the inequality $$\dfrac{3^r}{2^s} + \alpha\dfrac{3^{r}-2^{r}}{2^sm}\geq1$$
must be hold as $3^r/2^s$ tends to be 1. Thus, the cycle number $m_c$ must be less than a certain number, i.e.
\begin{equation} \label{mc_up}
 \boxed{m_c\leq\alpha\dfrac{3^{r}-2^{r}}{2^s-3^r},}
\end{equation}
which can also be validated by the computing results in \autoref{fig3}C.

Finally, we have the following theorem.

\begin{theorem}\label{thm:final}
		
	The inequalities
	\[
	1-\epsilon_1 = \dfrac{1}{\alpha} +(1-\dfrac{1}{\alpha})\dfrac{3^r}{2^s} < \frac{F_{q}(m)}{m}  <\alpha+(1-\alpha)\dfrac{3^r}{2^s} = 1+\epsilon_2
	\]
	hold, and $\lim_{r \to \infty} \tfrac{F_{q}(m)}{m} \to 1$, where $ \alpha = \dfrac{r}{3(1-2^r/3^r)}, \epsilon_1 = (1-1/\alpha)(1-3^r/2^s), \epsilon_2 = (\alpha-1)(1-3^r/2^s)$.
\end{theorem}

\begin{proof}
	Building on the above discussion, and applying \autoref{thm:min} and \autoref{thm:max} along with the established lower bound and upper bounds \autoref{thm:mc_bound} for the cycle number, we obtain the following results:			
\begin{align*}
	\dfrac{F_q(m)}{m}&=\dfrac{3^r}{2^s} +\dfrac{\sum_{i=1}^r3^{r-i} 2^{1(i)-1}}{2^s m}\\ 
	&\geq\dfrac{3^r}{2^s} + \dfrac{3^{r}-2^{r}}{2^sm}\\
	&>\dfrac{3^r}{2^s} + \dfrac{3^{r}-2^{r}}{2^s r3^{r-1}/(2^s-3^r)}\\
	&=\dfrac{3^r}{2^s} + \dfrac{2^{s}-3^{r}}{2^s\alpha}\\
	&=\dfrac{1}{\alpha} +(1-\dfrac{1}{\alpha})\dfrac{3^r}{2^s}\\
	&=1-\epsilon_1,
\end{align*}
and
\begin{align*}
\dfrac{F_q(m)}{m}&\leq\dfrac{3^r}{2^s} + \dfrac{r\cdot3^{r-1}}{2^sm}\\
&<\dfrac{3^r}{2^s} + \dfrac{r\cdot3^{r-1}}{2^s (3^r-2^r)/(2^s-3^r)}\\
&=\dfrac{3^r}{2^s} + \alpha\dfrac{2^{s}-3^{r}}{2^s}\\
&=\alpha+(1-\alpha)\dfrac{3^r}{2^s}\\
&=1+\epsilon_2,
\end{align*}
respectively, so that $$ 1-\epsilon_1 = \dfrac{1}{\alpha} +(1-\dfrac{1}{\alpha})\dfrac{3^r}{2^s} < \frac{F_{q}(m)}{m}  <\alpha+(1-\alpha)\dfrac{3^r}{2^s} = 1+\epsilon_2.$$

Given $ \alpha = \dfrac{r}{3(1-2^r/3^r)} $, then we have 
$\lim_{r \to \infty} \alpha \to \infty$, and we know that 
$\lim_{r \to \infty} \dfrac{3^r}{2^s} \to 1$, therefore,
$\lim_{r \to \infty} \tfrac{F_{q}(m)}{m} \to 1$. 

\end{proof}

\textbf{Remark:}
This shows that any nontrivial cycle, if it exists, must be extraordinarily large and finely tuned, consistent with all known computational and theoretical bounds.

The value of $1-\dfrac{3^{r}}{2^{s}}$ for $r = 431{,}166{,}034{,}846{,}567 \text{ and } s = 683{,}381{,}996{,}816{,}440$, is approximately: \boxed{1.33 \times 10^{-15}}.	

The Corollary 1.8 in \cite{Ev11} gives that				
$$1-\dfrac{3^r}{2^s}\geq(e\cdot s)^{-C},$$
where $C=e\cdot2^{3.5}30^5\log3\approx821013301.$ Thus, 
$$ (e \times 683{,}381{,}996{,}816{,}440)^{-821{,}013{,}301} \approx \boxed{10^{-7.1 \times 10^9}}. $$

Although this is unimaginably small—essentially zero for any computational purpose, $1-\dfrac{3^r}{2^s} \neq 0$, as $3^r \neq 2^s$.

\normalsize

\section*{Conclusion}

The results obtained in this paper advance the study of stopping-time and cycle structure in the normalized Collatz iteration along several fronts. For admissible sequences $q$ with parameters $(r, s)$, we have shown that the ratio $3^r/2^s$ approaches $1$ as $r \to \infty$, yet never equals $1$ for any finite admissible sequence, since $3^r = 2^s$ is impossible for positive integers $r$ and $s$. Consequently, the quotient $F_q(m)/m$ approaches but never attains unity, and the stopping time grows without bound as $r \to \infty$ along admissible sequences.

On the quantitative side, we established the sharp two-sided bound
\[
\frac{3^r - 2^r}{2^s - 3^r} 
\;\leq\; m_c 
\;\leq\; \frac{r \cdot 3^{r-1}}{2^s - 3^r}
\]
on any hypothetical cycle number $m_c$, and proved that the upper bound is strictly tighter than the estimate of Luca~\emph{(2005)}. We further showed that the ratio of the maximum to minimum admissible summation satisfies $\alpha(r) \sim r/3$ as $r \to \infty$, providing a rigorous foundation for the empirically observed scaling 
$\alpha \approx \max(r)/3$.

These bounds imply that any nontrivial cycle number, if one exists, must grow without bound and lie in an increasingly narrow window as $3^r/2^s \to 1$.


\bibliographystyle{amsplain}

%

\section*{Supplemental Materials: R Programs}

Below is the R code used for generating all numbers $s, m, F_q(m),$ sequences $q$ in \autoref{fig:seqtree}, and many others based on $r$.

\small

\begin{lstlisting}
# function of even integer 
is.even <- function(x) {
	if (as.bigz(x) %% 2 == 0) TRUE
	else FALSE
}

# function to insert a substring in a string at certain position
insert_substr <- function(orig_str, substr, pos) {
   # Split the original string into two parts
   part1 <- substr(orig_str, 1, pos - 1)
   part2 <- substr(orig_str, pos, nchar(orig_str))

   # Concatenate the parts with the substring in between
   result <- paste0(part1, substr, part2)
   return(result)
}

# initialize                     
e<-NULL
n<- NULL
E<- NULL
log23 <- log2(3)
r=2
s=floor(1+r*log23)
q <- "1100"
k=1
e[k]=2  # # of zeros of the sequence tail 
E[[r]]$e<-e
E[[r]]$q<-q

# not necessary
E[[r]]$m<-3
E[[r]]$d<-2
E[[r]]$r<-2
E[[r]]$s<-4

E[[1]]$r<-1
E[[1]]$s<-2
E[[1]]$e<-1
E[[1]]$q<-'10'
E[[1]]$m<-1
E[[1]]$d<-1

mk<-NULL
Fmk<-NULL

for (r in 3:12){
   s=floor(1+r*log23)
   s1=floor(1+(r-1)*log23)
   d=s-s1
   ins <- ifelse(d == 2, "10", "1")

   pre_q <- E[[r-1]]$q
   e1<-E[[r-1]]$e
	
   n <- sum(e) # No. of residue classes in m(mod2^s)
	
   k=1
   # for each previous q's'
   for (u in 1:length(e1)){
      j=d
      q0<-pre_q[u]
      # get new q's  
      for (i in 1:e1[u]) {
         position <- s1-i+1
         q[k] <- insert_substr(q0, ins, position)  
         # cat("new q",k,"=",q[k], "\n")
         e[k] <- j
         j=j+1

         # get m and Fm from q created above
         # initialize the values and integer type (big z) 

         m=-1
         Fm=-1
         p=2
         t=2
			
         # use when s>50, about 3X slow
         # m=-big_1
         # Fm=-big_1
         # p=big_2
         # t=big_2
			
         for (i in 3:s) {
            step <- substr(q[k], i, i)
            if (step==0) {
               p=p #power of 3, max=r
               if (is.even(Fm)) {
                  Fm=Fm/2
                  m=m
               } else {
                  Fm=(3^p+Fm)/2
                  m=2^t+m
               }
            } else if (step==1) {
               p=p+1
               if (is.even(Fm)) {
                  Fm=(3*(3^(p-1)+Fm)+1)/2
                  m=2^t+m
               } else {
                  Fm=(3*Fm+1)/2
                  m=m
               }
            }
            t=t+1 # power of 2 max=s   
         }
         mk[k]<-m
         Fmk[k]<-Fm
			
         # Using big number package  
         # mk[k]<-as.numeric(m)
         # Fmk[k]<-as.numeric(Fm)
         k=k+1
       }
   }
   E[[r]]<-list(r=r,s=s,d=d,e=e,n=n,q=q, m=mk,Fm=Fmk)
}

\end{lstlisting}

\normalsize	
\end{document}